\documentstyle[11pt]{article}
\input{pstricks.tex}
\input{mssymb.sty}
\newtheorem{thm}{Theorem}

\newtheorem{pro} {Proposition}
\newtheorem{dfn} {Definition}

\newenvironment{pf}{{\it Proof:}\quad}{\hfill$\Box$}

\begin{document}

\title{Guts of surfaces in punctured-torus bundles }
\author{Thilo Kuessner \\
Universit\"at Siegen, 
Mathematisches Institut, Flexstr.3, 57068 Siegen, \\
kuessner@mathuni-siegen.de}
\date{}
\maketitle
\begin{abstract}
\noindent
Let $M$ be a hyperbolic manifold of finite volume which fibers over the circle with fiber a once punctured 
torus, and let $S$ be an arbitrary incompressible surface in $M$. We determine the characteristic Jaco-Shalen-Johannson-submanifold
of $ M - S$ and show, in particular, that $Guts\left(M,S\right)$ is empty.
\end{abstract}

\section{Introduction}

Hyperbolic volume is an important homotopy invariant of hyperbolic manifolds and it is an interesting question to relate it to other topological invariants of (special classes of) hyperbolic 3-manifolds.

There is a topological definition of hyperbolic volume, due to Gromov(\cite{gro}). He defined a homotopy invariant $\parallel M\parallel$, the simplicial volume, and proved: $$\parallel M\parallel =\frac{1}{V_n}Vol\left(M\right)$$
if $M$ is an n-dimensional hyperbolic manifold of finite volume. ($V_n$ 
denotes the volume of a regular ideal straight n-simplex in hyperbolic n-space.)\\
There is an obvious inequality $$\parallel M\parallel\ge\frac{1}{n+1}\parallel\partial M\parallel$$ for any compact n-manifold, due to the fact that the boundary operator has norm $\le n+1$. If $M$ is a hyperbolic 3-manifold of finite volume, then $\parallel \partial M\parallel=0$, so this inequality gives no information.\\
Let $M$ be a hyperbolic manifold of finite volume and $S$ an arbitrary incompressible surface in $M$.
Agol used the simplicial volume to prove the inequality $$Vol\left(M\right)\ge - 2V_3 \chi\left(Guts\left(M,S\right)\right).$$
(We will explain the guts-terminology in the preliminaries, section 2.)
In particular, if $M_S:= M - {\mathcal{N}}\left(S\right)$ (that is, $M$ cut along $S$) admits a hyperbolic metric with totally geodesic boundary $\partial_1M_S:=S_1\cup S_2$ (where $S_1,S_2$ are the two copies of $S$ in $\partial M_S$)
and cusps corresponding to $\partial_0M_S:=\partial M_S - \partial_1M_S$, then
$$\parallel M\parallel\ge-\chi\left(S\right)=\frac{1}{4}\parallel\partial_1
M_S\parallel.$$\\
Agol used his inequality in \cite{ag2} to get useful estimates for the hyperbolic volume of 2-bridge links, and he suggested (in section 10) that the inequality should also be applied to other classes of 3-manifolds, where incompressible surfaces are
well-understood, namely alternating link complements and punctured-torus bundles.

For alternating link complements this task has been done in \cite{lac}. Lackenby showed the fairly surprising theorem that the hyperbolic 
volume of an alternating link complement is coarsely proportional to the twist number $t$ of the 
alternating diagram, precisely $\frac{1}{2}V_3\left(t-2\right)\le Vol\left(M\right)\le 10 V_3\left(t-1\right)$.
(The lower bound is exactly $-V_3\left(\chi\left(Guts\left(M,B\right)\right)
+\chi\left(Guts\left(M,W\right)\right)\right)$, where $B$ and 
$W$ are the checkerboard surfaces of a twist-reduced, prime, alternating link diagram. The upper bound is obtained by Agol and D.Thurston using ad-hoc methods.)

In this paper we determine the guts for an arbitrary incompressible surface in a hyperbolic once-punctured-torus bundle. 

\begin{thm}: Let $M$ be a once-punctured-torus bundle over ${\Bbb S}^1$ with 
hyperbolic monodromy. Let $S$ be an arbitrary connected, two-sided, properly embedded,
incompressible surface 
in $M$, not isotopic to the fiber nor the boundary torus. Associated to $\left(M,S\right)$ there is a natural number $k$ such that the following holds:\\
(i) if $k$ is even, then $k=-\chi\left(S\right) $ and the components of the characteristic submanifold of $M_S$ are $k$ regular neighborhoods of essential squares and the boundary annuli $\partial M_S\cap\partial M$, and $k$ Seifert fibered solid tori,\\
(ii) if $k$ is odd, then $2k=-\chi\left(S\right)$ and the components
of the characteristic submanifold of $M_S$ are $k$ regular neighborhoods of essential squares and the boundary annuli $\partial M_S\cap\partial M$, $k$
Seifert fibered solid tori, and one handlebody which is an $I$-bundle over
a one-sided, properly embedded surface.\\
In both cases, $Guts\left(M,S\right)$ is empty.
\end{thm}

The proof, which is given in section 3, proceeds by an elementary analysis
of the incompressible surfaces given by the classification due to Floyd-Hatcher
and Culler-Jaco-Rubinstein.

Theorem 1 is, of course, disappointing from the point of view of Agol's inequality. It says that Agol's inequality can not provide nontrivial estimates for the hyperbolic volume of once-punctured-torus bundles.

We mention that, on the other hand, we have an upper bound for 
the hyperbolic volume in terms of the 
monodromy. We denote by $t_l$ the Dehn-twist 
in the longitude and by $t_m$ the Dehn-twist in the meridian of $T$,
the once-punctured torus. It follows from Murasugi's theorem (\cite{cjr}, prop.1.3.1) that 
the mapping tori with monodromy $A$ and $B$ are homeomorphic if and only if 
$A$ is conjugate to $B^{\pm 1}$. 
Moreover, it is well-known that any mapping class
of the once punctured torus is conjugate to one of the form $A=\pm t_l^{l_1}t_m^{m_1}\ldots t_l^{l_r}t_m^{m_r}$, $l_1,\ldots,l_r,m_1\ldots,m_r\in{\Bbb N}$. (That one can get $l_1,\ldots,m_r\ge 1$ after suitable conjugation follows from the argument on page 267 of \cite{fh}.)
\begin{thm}: Let $M$ be a once-punctured-torus bundle over ${\Bbb S}^1$ with hyperbolic monodromy $A=
t_l^{l_1}
t_m^{m_1}\ldots t_l^{l_r}t_m^{m_r}$. Assume that $l_1,\ldots,l_r,m_1,\ldots,m_r\ge 1$. 
Then $$Vol\left(M\right)\le V_3\sum_{i=1}^r l_i+m_i.$$\end{thm}
\begin{pf} Let $A\in SL\left(2,{\Bbb Z}\right)$ be the monodromy of the bundle, and let $l$
%\gamma=\left\{\ldots,\frac{a_{-1}}{b_{-1}},\frac{a_0}{b_0},\frac{a_1}{b_1},\ldots\right\}$ 
be a minimal 
$A$-invariant path in the dual tree to the
Farey-triangulation. (All notions are explained in section 2.3.\ or in \cite{fh}.)
Associated to $l$ there is an ideal triangulation,
see the appendix
of \cite{fh}. Parker proved in \cite{par} that these ideal triangulation 
can be realized as a triangulation by straight ideal 
simplices (cf.\ also \cite{la1} for a stronger statement).
It is well-known that any straight ideal 3-simplex in hyperbolic 3-space 
has volume $\le V_3$ (cf.\ \cite{mi}). Hence $Vol\left(M\right)\le kV_3$,
for the number $k$ of tetrahedra in the ideal triangulation.

The tetrahedra in the ideal triangulation correspond precisely to the vertices of a minimal $A$-invariant path in the dual tree to the Farey triangulation. 
Each minimal $A$-invariant path belongs to the strip drawn in 
\cite{fh}, figure 3. It is clear from \cite{fh}, figure 3, that the number of vertices (modulo the action of $A$) is $l_1+m_1+\ldots +l_r+m_r$. This implies the claimed inequality.
\end{pf}\\

Equality in the inequality of Theorem 2 holds for powers of the Anosov cat map (i.e., finite covers of the figure-eight knot complement).

We emphasize that Theorem 2 is not true without the assumptions \\
$l_1,\ldots,l_r,m_1,\ldots,m_r\ge 1$.
\eject

\section{Preliminaries}
\subsection{3-manifolds with boundary-patterns}

We recall some terminology from \cite{joh}. 
Let $N$ be a compact $n$-manifold. A boundary-pattern for 
$N$ consists of a set $\underline{n}$ of compact, connected $\left(n-1\right)$-manifolds in $\partial N$, such that the intersection of any $i$ of them is a (possibly
empty) $\left(n-i\right)$-manifold.

A map $f:\left(M,\underline{m}\right)\rightarrow\left(N,\underline{n}\right)$
is called admissible if
$$\underline{m}=\bigcup_{G\in\underline{n}}\left\{\mbox{components of }f^{-1}G\right\}.$$
An admissible homotopy is an admissible map $f:\left(M\times I,\underline{m}\times I\right)\rightarrow\left(N,\underline{n}\right)$, where $\underline{m}\times I=\left\{G\times I: G\in\underline{m}\right\}$.

Let $\left(I,\underline{i}\right)$ be the closed intervall with the full boundary-pattern (i.e., consisting of both endpoints). 
An admissible curve in $N$ is an admissible mapping $k:\left(I,\underline{i}\right)
\rightarrow\left(N,\underline{n}\right)$ or $k:\left(S^1,\emptyset
\right)\rightarrow\left(N,\underline{n}\right)$. A $\partial$-compression disk for $k$ is 
an admissible mapping $g:\left(D,\underline{d}\right)\rightarrow\left(N,\underline{n}\right)$, where $D$ is a 2-disk with $j$ boundary faces, such that $k$ is one of the boundary faces and that the other $j-1$ boundary faces (apart from $k$) belong to $\underline{d}$, the boundary pattern of $D$. (The case $j=0$ is admitted if $k:\left(S^1,\emptyset\right)\rightarrow\left(N,\underline{n}\right)$. In this case, $k$ is the only boundary face of $D$.)The curve $k$ is called essential
if no compression disk exists. A submanifold $\left(F,\underline{f}\right)\subset \left(N,\underline{n}\right)$ is essential if essential curves in $F$ are also essential in $N$.

A boundary pattern $\underline{n}$ of a 3-manifold $N$ is 
called {\bf useful} if the boundary curve of any admissibly
embedded 2-disk $\left(D,\underline{d}\right)$,
where $\partial D$ has $i\le 3$ faces and $\partial D=\cup\left\{d: d\in\underline{d}\right\}$, bounds a disk $D^\prime$ in $\partial N$ such that $D^\prime\cap
\left(\cup_{G\in\underline{n}}\partial G\right)$ is the cone on $\partial D^\prime\cap
\left(\cup_{G\in\underline{n}}\partial G\right)$. 

In what follows, an {\bf annulus} is an annulus with the full boundary-pattern (i.e., consisting of two 
circles) and a {\bf square} is a 2-disk with four boundary faces such that the boundary-pattern consists of two 
non-adjacent boundary faces.

When talking about Seifert fibrations or $I$-bundles $p:W\rightarrow F$, we will always assume that the 
boundary pattern $\underline{w}$ of $W$ is such that 
$$\underline{w}=\left\{G; G=p^{-1}k \mbox{ for some }
k\in\underline{f} \mbox{ or }G\mbox{ a component of }\overline{\partial M - p^{-1}\partial F}\right\}$$ for some boundary pattern $\underline{f}$ of $F$.
An admissible {\bf F-manifold} in $\left(N,\underline{n}\right)$ is a submanifold $W\subset N$ such that each component of $W$ is an admissibly embedded $I$-bundle 
or Seifert fibre space and that $\overline{\partial W\setminus\partial N}$ 
is essential. (In particular, $\overline{\partial W\setminus\partial N}$ consists of 
essential squares, annuli and tori.) 

{\em Example:} Assume $M$ fibers over ${\Bbb S}^1$ with fiber $F$ and 
$\left(N,\underline{n}\right)=\left(M_F,F_1\cup F_2\right)$. Then $\left(N,\underline{n}\right)$ is an $I$-bundle over $F$. Indeed, $\underline{n}$
meets the above definition for the boundary pattern $\underline{f}=\emptyset$
of $F$. 
\eject

An admissible F-manifold $W\subset N$ 
is called {\bf characteristic} if for every component $V$ of $\overline{M\setminus W}$, $W\cup V$ 
is not an admissible F-manifold, and if every admissible F-manifold in $N$ can be admissibly isotoped into $W$.

{\em Remark:} Assume that the complement $\partial_0N$ of the boundary pattern 
in $\partial N$ consists of annuli. The regular neighborhood of 
any essential square 
$D$ can be considered as an I-bundle (respecting boundary pattern)
by projecting to the product neighborhood of one of those faces of $D$
which belong to the boundary pattern of $D$. (This neighborhood is the base surface of the $I$-bundle, it has empty boundary pattern.) We remark that this is not a maximal $I$-bundle. 
Indeed, the 
$I$-bundle projection can be extended (respecting boundary pattern) in the obvious way to a product neighborhood of those annuli in $\partial_0N$ which intersect $D$.
These larger $I$-bundles have annuli (not squares) as their boundary 
components. Since this can be done for any essential square, one sees that
the boundary of the characteristic submanifold consists of essential tori and annuli only.

The following is Proposition 9.4.\ and Corollary 10.9.\ in \cite{joh}, see also \cite{js}.
\begin{pro}(Jaco-Shalen, Johannson): Let $\left(N,\underline{n}\right)$ be an irreducible, boundary-irreducible 3-manifold with 
useful boundary pattern. Then there exists a characteristic submanifold in $\left(N,\underline{n}\right)$, unique up to admissible isotopy.\end{pro}

\subsection{Guts of surfaces}

Let $M$ be a compact 3-manifold (with boundary) and $S$ a two-sided properly embedded surface in $M$.
Let $N=M_S:=\overline{M-{\mathcal{N}}\left(S\right)}$ be the manifold obtained by 
splitting $M$ along $S$, that is, the closure of the complement of a
regular neighborhood ${\mathcal{N}}\left(S\right)$ of $S$ in $M$.

$\partial N$ contains two copies of $S$, which we denote $S_1$ and $S_2$. We consider the boundary pattern $\underline{n}=S_1\cup S_2$. If $S$ was incompressible, then $\underline{n}$ is obviously a useful boundary-pattern.
\begin{dfn} Let $M$ be a compact 3-manifold (with boundary) and $S$ a two-sided, properly embedded
surface in $M$. Assume that $\left(N,\underline{n}\right):=\left(M_S,S_1\cup S_2\right)$ satisfies the assumptions of Proposition 1. Then we define $Guts\left(M,S\right)$ to be the closure of the complement
of the characteristic submanifold of $\left(N,\underline{n}\right)$ in $N$.
\end{dfn}
 
If $M$ is irreducible and boundary-irreducible, and $S$ 
incompressible, then it can be deduced from Thurston's geometrization 
theorem for Haken manifolds that $Guts\left(M,S\right)$ admits a
hyperbolic metric with $\partial Guts\left(M,S\right) \cap \left(S_1\cup S_2\right)$ totally geodesic.

\subsection{Incompressible surfaces in hyperbolic once-punctured-torus bundles}

%We denote $T^2$ the 2-torus and $T:=T^2-D^2$ the complement of an embedded disk in $T^2$. ($T$ is called a once punctured torus.) Elements of $SL\left(2,{\mathbb Z}\right)$ act linearly on $T^2={\mathb R}^2/{\mathbb Z}^2$, preserving $0$. Thus we have a well-defined action of $SL\left(2,{\mathbbZ}\right)$ on 
It is well-known that the mapping class group of the once punctured torus is $SL\left(2,{\Bbb Z}\right)$, representatives acting linearly on $T^2-0=\left({\Bbb R}^2-{\Bbb Z}^2\right)/{\Bbb Z}^2$. Hence,
a once-punctured-torus bundle $M$ is determined by its monodromy $A\in SL\left(2,{\Bbb Z}\right)$. $M$ is hyperbolic if and only if $A$ is hyperbolic (i.e., diagonalizable and not $\pm\ identity$).

\begin{dfn}: Let $T=T^2-D^2$ be
the once punctured torus and $I=\left[a,b\right]$ a closed intervall.
A {\bf standard saddle} of bottom slope $\frac{p}{q}$ and top slope $\frac{r}{s}$ in $T\times I$ is an 
embedded 8-gon in $T\times I$ with faces (in consecutive order) $e_1,\ldots,e_8$ such that \\
- $e_1$ and $e_5$ are arcs of slope $\frac{p}{q}$ in $T\times\left\{a\right\}$,\\
- $e_3$ and $e_7$ are arcs of slope $\frac{r}{s}$ in $T\times\left\{b\right\}$,\\
- for $i=2,4,6,8$ is $e_i$ of the form $\left\{x_i\right\}\times I$ for some $x_i\in\partial T$.\end{dfn}

Incompressible surfaces in once-punctured-torus bundles have been classified
in \cite{cjr} and (for hyperbolic $M$) in \cite{fh}.      
We use the approach of \cite{fh}.

Consider the ideal Farey-triangulation of the upper half-plane model of the hyperbolic plane ${\Bbb H}^2$, which is as follows:
the ideal vertices are ${\Bbb Q}$ and $\infty=\frac{1}{0}$, two ideal vertices $\frac{a}{b}$ and $\frac{c}{d}$ are connected by an edge if $ad-bc=\pm 1$.

$SL\left(2,{\Bbb Z}\right)$ acts on ${\Bbb Q}\cup\left\{\infty\right\}$ 
by fractional linear transformations, this action extends to an action 
preserving the Farey-triangulation. 
Let $A\in SL\left(2,{\Bbb Z}\right)$. Consider a minimal $A$-invariant edge-path $\gamma$ in the Farey-triangulation. 
Let $\ldots, \frac{a_{-1}}{b_{-1}},\frac{a_0}{b_0},\frac{a_1}{b_1},\ldots$ be the vertices of $\gamma$. Since $\gamma$ is $A$-invariant, we have some $k\in{\Bbb N}$
such that $A\left(\frac{a_i}{b_i}\right)=\frac{a_{i+k}}{b_{i+k}}$ for all $i$.

Let $M$ be the once-punctured-torus bundle with monodromy $A$. To each minimal $A$-invariant edge-path $\gamma$ we associate an incompressible surface $S_\gamma\subset M$ as follows.

Consider $M$ as the quotient of $T\times\left[0,1\right]$ under the identification $\left(x,0\right)\sim\left(A\left(x\right),1\right)$. We consider the surface in $T\times\left[0,1\right]$ whose intersection with $T\times\left[\frac{i}{k},
\frac{i+1}{k}\right]$ is the standard saddle with bottom slope $\frac{a_i}{b_i}$ 
and top slope $\frac{a_{i+1}}{b_{i+1}}$, for $i=0,1,\ldots,k-1$. Since $A\left(\frac{a_0}{b_0}\right)=\frac{a_k}{b_k}$, this gives a properly embedded surface in $M$.
By \cite{fh}, $S$ is incompressible in $M$. If $k$ is even, then $S$ is two-sided
and we let $S_\gamma:=S$.
If $k$ is odd, then $S$ is one-sided and we let $S_\gamma$ be the boundary of a
regular neighborhood of $S$. 

The main result of \cite{fh} is that each incompressible two-sided surface in $M$ (apart from the fiber and the boundary torus)
arises this way.

\section {Proof of Theorem 1}

Let $M$ be the once-punctured-torus bundle with hyperbolic monodromy $A$. Let $S_\gamma$ be the incompressible two-sided surface associated 
to the $A$-invariant minimal edge-path $\gamma$, as explained in section 
2.3. We distinguish the cases that $k$, the period of $\gamma$, is even or odd. \\

{\bf Case 1}: $k$ even.\\
That is, if $\gamma$ has vertices $\left\{\frac{a_i}{b_i}:i\in{\Bbb Z}\right\}$ with 
$A\left(\frac{a_i}{b_i}\right)=\frac{a_{i+k}}{b_{i+k}}$ for a fixed $k\in{\Bbb N}$ and all $i\in{\Bbb Z}$, 
then $S_\gamma$ is composed of the $k$ standard saddles with bottom slopes $\frac{a_i}{b_i}$
and top slopes $\frac{a_{i+1}}{b_{i+1}}$, for $i=0,1,\ldots,k-1$. Note that $\chi\left(S_\gamma\right)=-k$.

Let $N=M_S$. A family of squares in $N$ is constructed as follows. Consider $M$ as a quotient of $T\times\left[0,1\right]$ in the obvious way.
For each $i=0,1,\ldots,k-1$ let $D_i$ be the square embedded 
in $T\times\left\{\frac{i}{k}\right\}$ such that two non-adjacent boundary faces are arcs 
of slope $\frac{a_i}{b_i}$ and the two other non-adjacent boundary faces belong to $\partial T\times\left\{\frac{i}{k}\right\}$. These squares in $M$ give squares in the manifold $M_S$ with boundary-pattern $S_1\cup S_2$ (see section 2.2.). 

Note that the squares $D_1,\ldots,D_k$ are essential, because the boundary faces not belonging to the 
boundary-pattern are essential in $T$ and thus (since $T$ is essential in $M$) also in $M$ and $M_S$. Thus $D_1\cup\ldots\cup D_k$ 
can be isotoped into the characteristic submanifold.
Let $V$ be a regular neighborhood of $D_1\cup\ldots\cup D_k$.
(We may, as described in the remark in section 2.1., enlargen each component of $V$ as an I-bundle by a neighborhood of the boundary annulus. This does not change 
neither the homeomorphism type of the complement nor its boundary
pattern.)\\

{\bf Claim:} {\em The complement of $V$ consists of $k$ solid tori.}\\
\\
This should be fairly clear from the picture of the standard saddles. To make a precise argument,
consider the once-punctured torus $T$ as pictured below, with the pairs of $m$'s resp.\ $l$'s to
be identified. The standard saddle $K$ with bottom slope $\frac{0}{1}$ and top slope $\frac{1}{0}$
is embedded in $T\times \left[\frac{i}{k},\frac{i+1}{k}\right]
$ as the 8-gon with edges $l_-\times \left\{\frac{i}{k}\right\},\left\{p_1\right\}\times \left[
\frac{i}{k},\frac{i+1}{k}\right], m_+\times \left\{\frac{i+1}{k}\right\},
\left\{p_2\right\}\times \left[
\frac{i}{k},\frac{i+1}{k}\right], l_+\times\left\{\frac{i}{k}\right\},
\left\{p_3\right\}\times \left[
\frac{i}{k},\frac{i+1}{k}\right],
m_-\times\left\{\frac{i+1}{k}\right\},
\left\{p_4\right\}\times \left[
\frac{i}{k},\frac{i+1}{k}\right].$ Let $D$ be the disk, 
containing $m\times\left\{\frac{i+1}{k}\right\}$, which has $m_-\times
\left\{\frac{i+1}{k}\right\}$ and $m_+\times
\left\{\frac{i+1}{k}\right\}$ as two of its boundary faces, and its other two boundary faces on $\partial T\times
\left\{\frac{i+1}{k}\right\}$. Then $K\cup D$ decomposes $T\times \left[\frac{i}{k},
\frac{i+1}{k}\right]$ into two connected components. Both are solid tori, indeed one of them is obtained from
a cube by identifying two copies of $m\times \left[\frac{i}{k},
\frac{i+1}{k}\right]$, the other is obtained from a
cube by identifying two copies of $l\times
\left[\frac{i}{k},
\frac{i+1}{k}\right]$. 

Now, apply $$f:=\left(\begin{array}{cc}  a_i &   {a_{i+1}} \\ b_i &   b_{i+1}\end{array}\right)\in SL\left(2,{\Bbb Z}\right)$$ to
(the first factor of) $T\times\left[\frac{i}{k},\frac{i+1}{k}\right]$. Then $K$ is mapped to 
the standard saddle of bottom slope $\frac{b_i}{a_i}$ and top slope $\frac{b_{i+1}}{a_{i+1}}$
and $D$ is mapped to $D_{i+1}$. Since $f$ is a homeomorphism, the images of the components of $T\times\left[
\frac{i}{k},\frac{i+1}{k}\right]-\left(K\cup D\right)$ are solid tori, i.e., $D_{i+1}$ decomposes the complement of
the standard saddle of boundary slopes $\frac{b_i}{a_i}$ and $\frac{b_{i+1}}{a_{i+1}}$ in 
$T\times\left[
\frac{i}{k},\frac{i+1}{k}\right]$ into two solid tori. So, if we consider one of the connected components of 
$M_S - V$, e.g.\ the component bounded by $D_i$ and $D_{i+1}$, i.e.\ intersecting $T\times\left[
\frac{i}{k},\frac{i+1}{k}\right]$ and $T\times\left[
\frac{i-1}{k},\frac{i}{k}\right]$, then this component is obtained from two solid tori (the intersections with
$T\times\left[
\frac{i}{k},\frac{i+1}{k}\right]$ resp.\ $T\times\left[
\frac{i-1}{k},\frac{i}{k}\right]$) by glueing them along a
al annulus (which is the complement in $T\times\left\{\frac{i}{k}\right\}$ of the disk containing $m\times\left\{\frac{i}{k}\right\}$ bounded by 
$m_+\times\left\{\frac{i}{k}\right\}$, two arcs in $\partial T\times\left\{\frac{i}{k}\right\}$,
and $m_-\times\left\{\frac{i}{k}\right\}$.) Glueing the
two solid tori along 
the longitudinal annulus yields a new solid torus. 

It remains to show that this solid torus is Seifert fibered, respecting the boundary pattern.
First remark that, for a solid torus with boundary pattern an annulus, a Seifert fibration in which this boundary pattern consists of boundary fibers exists if and only if this annulus is not meridional,
i.e., winds along the core at least once. In our case, the intersection of the boundary pattern
with $T\times\left[\frac{i}{k},\frac{i+1}{k}\right]$ is half of a
annulus of slope $\frac{b_{i+1}}{a_{i+1}}$, and the intersection of the boundary pattern
with $T\times\left[\frac{i-1}{k},\frac{i}{k}\right]$ is half of a
annulus of slope $\frac{b_i}{a_i}$. These two half-annuli are glued together to
form an annulus. Recall that the core of the solid torus has slope $\frac{b_i}{a_i}$.
Hence the annulus in the boundary pattern is not meridional and there exists a Seifert fibration.

This shows the claim if $k$ was even.\\

{\bf Case 2:} $k$ is odd. \\
In this case, $S$ is constructed as above, but $S_\gamma$ is the (connected) boundary of a standard neighborhood 
${\mathcal{N}}\left(S\right)$ of $S$. Note that $\chi\left(S_\gamma\right)=-2k$.\\
The same argument as in case 1 shows that the complement of 
${\mathcal{N}}\left(S\right)$ can be cut 
along $k$ essential squares into $k$ solid tori. So it remains
to look at ${\mathcal{N}}\left(S\right) - {\mathcal{N}}\left(S_\gamma\right)={\mathcal{N}}\left(S\right)
 - {\mathcal{N}}\left(\partial {\mathcal{N}}\left(S
\right)\right)$. Now, $S$ is a nonorientable surface of Euler characteristic $-k$ with one
boundary component, thus ${\mathcal{N}}\left(S\right) - {\mathcal{N}}\left(
\partial{\mathcal{N}}\left(S\right)\right)$ is a 
handlebody which is an $I$-bundle over $S$.
\hfill$\Box$\\
\\
\\
\\
\\
\setlength{\unitlength}{4mm}
\begin{picture}(8,5.5)(-8,-4)
\put(-4,4){\line(1,0){8}}
\put(-4,-4){\line(1,0){8}}
\put(-8,-2){\line(0,1){4}}
\put(8,-2){\line(0,1){4}}
\put(-6,3){\line(1,0){12}}
\put(-6,-3){\line(1,0){12}}
\put(-6,-3){\line(0,1){6}}
\put(6,-3){\line(0,1){6}}
\put(-8,2){\line(2,1){4}}
\put(4,-4){\line(2,1){4}}
\put(-4,-4){\line(-2,1){4}}
\put(8,2){\line(-2,1){4}}

\put(8.3,0){$m$}
\put(-7.8,0){$m$}
\put(0,3.2){$l$}
\put(0,-4.8){$l$}
\put(6.3,0){$m_+$}
\put(-5.8,0){$m_-$}
\put(0,2.2){$l_+$}
\put(0,-3.8){$l_-$}

\put(6.4,3.4){$p_2$}
\put(6.4,-3.4){$p_1$}
\put(-6.4,-3.7){$p_4$}
\put(-6.4,3.4){$p_3$}

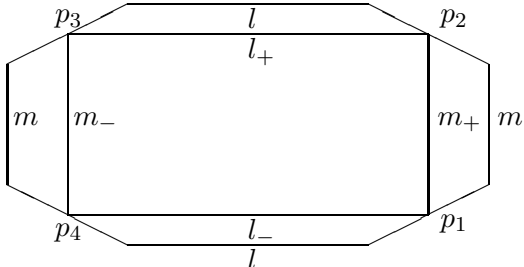
\put(-4.4,-5.7){Pairwise identification of $m$'s and $l$'s yields punctured torus} 
\end{picture}\\
\\
\\
\\
{\bf Question:}
It seems somewhat astonishing to me that in all cases considered so far (2-bridge links,
alternating links, once-punctured-torus bundles), the characteristic subpairs
of $\left(M,S\right)$ consist only of regular neighborhoods of 
essential squares (and boundary annuli)
and of Seifert fibered solid tori. 
The squares arise from apparently different reasons: 
If $S$ is the spanning surface of a two-bridge link $K_{\frac{p}{q}}$, then the squares correspond
to the '2'-s in the associated continued fraction expansion, i.e., to the occurences
of $b_i=2$ in the expansion $\frac{p}{q}=r+\frac{1}{b_1-\frac{1}{b_2-\frac{1}{b_3-\ldots}}}$. If $S$ is the black checkerboard surface for an alternating link complement, then the squares correspond (not one-to-one)
to the white bigon regions in the alternating diagram. If $S$ 
is an arbitrary incompressible surface in a once-punctured-torus bundle, then the squares correspond to the vertices of the associated  minimal invariant edge-path in the Farey-triangulation.\\
This naturally raises the question, how general
this phenomenon may be among arbitrary hyperbolic 3-manifolds, i.e. whether there is some unified reason restricting the existence of essential annuli in $M-S$ and explaining the occurence of squares.\\
\\
{\bf Higher genus.} The analogue of corollary 1 for surface bundles of fiber 
genus $\ge 2$ is not true. Counterexamples can be found in \cite{ht}. Namely, 
it follows from \cite{ht} (the remark after the corollary on p.228), 
that, for example, the complement of the knot $K_{\frac{4}{11}}$ 
fibers over the circle. On the other hand, \cite{ht} constructs incompressible surfaces $S$ corresponding 
to the continued fraction expansion $\frac{4}{11}=\left[3,4\right]$ and then $\chi\left(Guts\left(M,S\right)\right)=-1$ follows from the computation in
\cite{ag2}, section 7. This negative result is, of course, positive from the
point of view of hyperbolic volume. It says that Agol's inequality provides nontrivial estimates for surface bundles of higher fiber genus. However, to explore this in a systematic way, one would need a classification of incompressible surfaces in surface bundles.


\begin{thebibliography}{18}

\bibitem{ag1} AGOL (I.). \textsc{Topology of hyperbolic 3-manifolds}, PhD-Thesis UCSD, (1998), http://www.math.uic.edu/~agol/shortthesis.ps
\bibitem{ag2} AGOL (I.). \textsc{Lower bounds on volumes of hyperbolic Haken manifolds}, Preprint, http://front.math.ucdavis.edu/math.GT/9906182.

\bibitem{cjr}CULLER (M.), JACO (W.) and RUBINSTEIN (H.). - \textsc{Incompressible surfaces in once-punctured torus bundles}, Proc.\ London Math.\ Soc., {\bf 45}, (1982), pp.385-419.


\bibitem{fh}FLOYD (W.) and HATCHER (A.). - \textsc{Incompressible surfaces in punctured-torus bundles}, Top.\ Appl., {\bf 13}, (1982), pp.263-282.
\bibitem{gro}GROMOV (M.). - \textsc{Volume and bounded cohomology}, Publ.\ IHES, {\bf 56}, (1982), pp.5-100.
\bibitem{ht}HATCHER (A.) and THURSTON (W.). - \textsc{Incompressible surfaces in 2-bridge knot complements}, Invent.\ Math., {\bf 79}, (1985), pp.225-246.
\bibitem{js}JACO (W.) and SHALEN (P.). - \textsc{Seifert fibered spaces in 3-manifolds}, Mem.\ Amer.\ Math.\ Soc., {\bf 21}, (1979), no.220.
\bibitem{joh}JOHANNSON (K.). - \textsc{Homotopy Equivalences of 3-Manifolds with Boundaries}, Lecture Notes in Mathematics 761, Springer Verlag, (1979).
\bibitem{la1}LACKENBY (M.). - \textsc{The canonical decomposition of once-punctured torus bundles}, Comm.\ Math.\ Helv., {\bf 78}, (2003), pp.363-384.
\bibitem{lac}LACKENBY (M.). - \textsc{The volume of hyperbolic alternating link complements}, Proc.\ London Math.\ Soc., {\bf 88}, (2004), pp.204-224.       
\bibitem{mi}MILNOR (J.). - \textsc{Computation of volume }, http://www.msri.org/publications/books/gt3m/PS/7.ps
\bibitem{ns}NEUMANN (W.) and SWARUP (G.). - \textsc{Canonical decompositions of 3-manifolds}, Geom.\ Topol., {\bf 1}, (1997), pp.21-40.
\bibitem{par}PARKER (J.). - \textsc{Tetrahedral decompositions of punctured torus bundles}, London Math.\ Soc.\ Lecture Notes, {\bf 299}, (2003), pp.275-291.
\end{thebibliography}
\end{document}